\documentclass[12pt,english]{article}
\usepackage[T1]{fontenc}
\usepackage[latin9]{inputenc}
\usepackage[a4paper]{geometry}
\usepackage{color}
\usepackage{float}
\usepackage{amsthm}
\usepackage{amsmath}
\usepackage{amssymb}
\usepackage{graphicx}

\makeatletter

\floatstyle{ruled}
\newfloat{algorithm}{tbp}{loa}
\providecommand{\algorithmname}{Algorithm}
\floatname{algorithm}{\protect\algorithmname}

\theoremstyle{plain}
\newtheorem{thm}{\protect\theoremname}
\theoremstyle{plain}
\newtheorem{cor}[thm]{\protect\corollaryname}
\theoremstyle{definition}
\newtheorem{defn}[thm]{\protect\definitionname}
\theoremstyle{plain}
\newtheorem{lem}[thm]{\protect\lemmaname}

\usepackage{algorithmic}
\usepackage{amsfonts}
\usepackage{url}

\makeatother

\usepackage{babel}
\providecommand{\corollaryname}{Corollary}
\providecommand{\definitionname}{Definition}
\providecommand{\lemmaname}{Lemma}
\providecommand{\theoremname}{Theorem}

\begin{document}

\title{High-Performance Algorithms for Computing the Sign Function of Triangular
Matrices}

\author{Vadim Stotland, Oded Schwartz, and Sivan Toledo}
\maketitle
\begin{abstract}
Algorithms and implementations for computing the sign function of
a triangular matrix are fundamental building blocks in algorithms
for computing the sign of arbitrary square real or complex matrices.
We present novel recursive and cache efficient algorithms that are
based on Higham's stabilized specialization of Parlett's substitution
algorithm for computing the sign of a triangular matrix. We show that
the new recursive algorithms are asymptotically optimal in terms of
the number of cache misses that they generate. One of the novel algorithms
that we present performs more arithmetic than the non-recursive version,
but this allows it to benefit from calling highly-optimized matrix-multiplication
routines; the other performs the same number of operations as the
non-recursive version, but it uses custom computational kernels instead.
We present implementations of both, as well as a cache-efficient implementation
of a block version of Parlett's algorithm. Our experiments show that
the blocked and recursive versions are much faster than the previous
algorithms, and that the inertia strongly influences their relative
performance, as predicted by our analysis.
\end{abstract}

\section{\label{sec:Introduction}Introduction}

The sign of a square complex matrix $A$ is defined by extending the
scalar function
\[
\text{sign}(z)=\text{sign}(x+iy)=\left\{ \begin{array}{ll}
1 & x>0\\
-1\quad & x<0.
\end{array}\right.
\]
to matrices. For a diagonalizable matrix $A=ZDZ^{-1}$ the sign can
be defined by applying $\text{sign}(z)$ to the eigenvalues of $A$,
\[
\text{sign}(A)=Z\begin{bmatrix}\text{sign}(d_{11})\\
 & \text{sign}(d_{22})\\
 &  & \ddots\\
 &  &  & \text{sign}(d_{nn})
\end{bmatrix}Z^{-1}\;;
\]
the definition can be extended to the non-diagonalizable case in a
variety of equivalent ways~\cite[Section~1.2]{HighamFoM}; from here
on, we use the term \emph{function} to refer to a mapping that satisfies
these equivalent definitions. The matrix sign function is not defined
when $A$ has purely imaginary eigenvalues (and is clearly ill-conditioned
on matrices with eigenvalues that are almost imaginary).

One way to compute the sign function is to first compute a Schur decomposition
of $A=QTQ^{*}$, where $T$ is upper triangular and $Q$ is unitary,
then compute $U=\text{sign}(T)$, and finally form $\text{sign}(A)=QUQ^{*}$.
In this paper we focus on computing the sign function of a triangular
matrix, which can be used a a building block in an algorithm for general
matrices. Parlett discovered a substitution-type algorithm that can
compute many functions of triangualr matrices~\cite{Parlett:1976:RAE}.
The algorithm exploits the equation $UT=TU$ that for any function
$U$ of $T$ satisfies and the fact that if $T$ is triangular, so
is $U$. Parlett's technique breaks down when $T$ has repeated eigenvalues
(and becomes unstable when it has clustered eigenvalues). 

Higham proposed an improved version that we refer to as the \emph{Parlett-Higham}
algorithm, which applies only for the sign function, and which avoids
breakdowns~\cite[Algorithm~5.5]{HighamFoM}. A more generic way to
avoid breakdowns and instability in Parlett's algorithm is to reorder
the Schur form so that eigenvalues are clustered along the diagonal
of $T$ and to apply a block version of Parlett's substitution~\cite{Davies:2003:SPA}.
This approach requires some other way to compute the sign of diagonal
blocks of $T$; the off-diagonal blocks are computed by solving Sylvester
equations. We refer to this method as the \emph{Parlett-Sylvester}
technique. The algorithm that computes the sign of diagonal blocks
must be able to cope with a clustered spectrum (up to the case of
repeated eigenvalues); Parlett's method cannot usually be applied
to these blocks. However, in the case of the sign function, clustering
the eigenvalues according to their sign provides a trivial way to
construct the two diagonal blocks of $U$: one is identity $I$ and
the other a negated identity $-I$, usually of different dimensions.

\paragraph*{Our Contributions}

This paper presents high-performance algorithms for computing the
sign of a triangular matrix. To obtain high performance, we take two
measures. First, we choose whether to use the Parlett-Higham substitution
algorithm or the Parlett-Sylvester algorithm, by estimating the amount
of work each of them would require. We show their complexity may differ
asymptotically, hence choosing the right one is essential. Second,
we reorder the operations that our algorithms perform, so as to reduce
cache misses and inter-processor communication. The reordering techniques
apply both to the Parlett-Higham and to the Parlett-Sylvester algorithms.

\paragraph*{Paper Organization}

The rest of the paper is organized as follows. Section~\ref{sec:Background}
presents the basic Parlett recurrence for functions of triangular
matrices as well as Higham's stabilized version for the sign function
and the Parlett-Sylvester approach. Section~\ref{sec:Arithmetic-Efficiency}
analyzes the number of arithmetic operations that the two approaches
perform and show that the Parlett-Sylvester is less efficient when
the inertia is balanced but much more efficient when it is not. Section~\ref{sec:Communication-Lower-Bounds}
presents lower bounds on the asymptotic number of cache misses that
these algorithms much generate. Section~\ref{sec:Communication-Efficient-Algs}
presents recursive cache-efficient variants of the Parlett-Higham
algorithm, which are asymptotically optimal by the previous section.
Section~\ref{sec:Experimental-Results} shows that the new algorithms
and our implementation of the Parlett-Sylvester algorithm are indeed
fast and that their performance in practice matches our theoretical
predictions. We presents our conclusions in Section~\ref{sec:Conclusions}.

\section{\label{sec:Background}Background}

Any matrix function $F=\phi(A)$ commutes with its argument, $AF=FA$.
The function $F=\phi(T)$ of an upper triangular matrix $T$ is also
triangular. Parlett used these facts to construct a substitution-type
algorithm to compute $F=\phi(T)$. By rearranging the expression for
the $i,j$ element in the product $TF=FT$
\[
\sum_{k=i}^{j}t_{ik}f_{kj}=\sum_{k=i}^{j}f_{ik}t_{kj}\;,
\]
where $t_{ik}$ is the $i,k$ element of $k$ and so on, we can almost
isolate $f_{ij}$
\begin{equation}
(t_{ii}-t_{jj})f_{ij}=f_{ii}t_{ij}-f_{jj}t_{ij}+\sum_{k=i+1}^{j-1}\left(f_{ik}t_{kj}-t_{ik}f_{kj}\right)\;.\label{eq:element-parlett}
\end{equation}
This allows us to obtain the value of $f_{ij}$ as a function of $f_{ik}$
for $k<j$ and $f_{kj}$ for $k<j$. These equations do not constrain
the diagonal elements of $F$ (the equations are $t_{ii}f_{ii}=f_{ii}t_{ii}$),
but it is easy to see that they must satisfy $f_{ii}=\phi(t_{ii})$.
The complete algorithm is shown in Algorithm~\ref{alg:parlett}.

\begin{algorithm}
\protect\caption{\label{alg:parlett}Parlett's substitution algorithm to compute a
function of a lower triangular matrix $T\in\mathbb{C}^{n\times n}$
with distinct diagonal elements (eigenvalues). }

\begin{algorithmic}[1]

\STATE\textbf{for} \textbf{$i=1:n,\; f_{ii}=\phi(t_{ii})$}

\STATE\textbf{for} $j=2:n$

\STATE\textbf{~~~for $i=j-1:-1:1$}

\STATE~~~~~~$f_{ij}=\frac{1}{t_{ii}-t_{jj}}\left(t_{ij}\left(f_{ii}-f_{jj}\right)+\left(\sum_{k=i+1}^{j-1}f_{ik}t_{kj}-t_{ik}f_{kj}\right)\right)$

\STATE\textbf{~~~end}

\STATE\textbf{end}

\end{algorithmic}
\end{algorithm}

Clearly, the algorithm breaks down if $T$ has repeated eigenavlues
($t_{ii}=t_{jj}$ for some $i$ and $j$). Pairs of nearby but unequal
eigenvalues (small $|t_{ii}-t_{jj}|$) tend to cause growth in $F$
because of divisions by small quantities. In some cases this is related
to ill conditioning of $F$, but not always. In some cases, the growth
is associated with an instability in the algorithm rather than with
poor conditioning. 

One way to address this issue, at least partially, is to partition
$F$ and $T$ into blocks and write the corresponding block-matrix-multiplication
equations that $TF=FT$ defines~(\cite{Parlett:1974:CFT}, cited
by~\cite{HighamFoM}). The partitioning is into square diagonal blocks
and possibly rectangular off-diagonal blocks. In this version we cannot
isolate off-diagonal blocks $F_{ij}$ because they do not necessarily
commute with diagonal blocks of $T$, so the equations that define
$F_{ij}$ are not simple substitution-type equations but rather they
form Sylvester equations, as shown in Algorithm~\ref{alg:parlett-sylvester}. 

\begin{algorithm}
\protect\caption{\label{alg:parlett-sylvester}The Parlett-Sylvester substitution algorithm
to compute a matrix function given a partitioning of the row and column
indices into $m$ blocks.}

\begin{algorithmic}[1]

\STATE\textbf{for }$i=1:m,\; F_{ii}=\phi(T_{ii})$

\STATE\textbf{for }$j=2:m$

\STATE~~~\textbf{for }$i=j-1:-1:1$

\STATE~~~~~~Solve for $F_{ij}$the Sylvester equation

~~~~~~$T_{ii}F_{ij}-F_{ij}T_{jj}=F_{ii}T_{ij}-T_{ij}F_{jj}+\sum_{k=i+1}^{j-1}\left(F_{ik}T_{kj}-T_{ik}F_{kj}\right)$

\STATE~~~\textbf{end}

\STATE\textbf{end}

\end{algorithmic}
\end{algorithm}

Here too, the equations that drive the algorithm say nothing about
diagonal blocks $F_{ii}$, so they must be computed in some other
way; we discuss this later. The Sylvester equation for $F_{ij}$ is
singular if $T_{ii}$ and $T_{jj}$ have common eigenvalues, and are
ill-conditioned if they have nearby eigenvalues. Hence, for this method
to work well, the partitioning of $F$ and $T$ needs to be such that
different diagonal blocks of $T$ share no common eigenvalues, and
ideally, do not have nearby eigenvalues.

Davies and Higham proposed a framework that uses this approach for
essentially any function $\phi$~\cite{Davies:2003:SPA}. Their framework
begins by clustering of eigenvalues of $T$. The clusters are made
as small as possible under the condition that they are well separated.
Note that if eigenvalues are highly clustered, the framework may end
up with a single large cluster. This is undesirable from the computational
complexity perspective, but avoids numerical problems. The framework
then uses an algorithm by Bai and Demmel~\cite{Bai:1993:SDB} to
reorder $T$ unitarily so as to make the eigenvalues in each cluster
adjacent, $T=Q\tilde{T}Q^{*}$. The Parlett-Sylvester algorithm then
computes the function $\tilde{F}=\phi(\tilde{T})$, which is transformed
back into the function $F$ of $T$, $F=Q\tilde{F}Q^{*}$. The diagonal
blocks of $\tilde{F}$ cannot be computed by a Parlett recurrence
because the diagonal blocks of $\tilde{T}$ have clustered or repeated
eigenvalues. Davies and Higham proposed that a Pade approximation
be used to compute these blocks. The Pade approach is very general
but becomes very expensive if diagonal blocks are large. 

However, in the special case of the sign function we can partition
the eigenvalues by the sign of their real part. In this case, the
functions of the two resulting diagonal blocks of $\tilde{T}$ are
trivial: the identity is the sign of the block with the positive eigenvalues
(right half of the complex plane) and a negated identity is the sign
of the block with the negative eigenvalues (left half of the plane)~\cite[Section~5.2]{HighamFoM}.

Higham also proposed another specialization of Parlett's method to
the sign function~\cite[Algorithm~5.5]{HighamFoM}. The matrix sign
$U=\text{sign}(T)$ satisfies another matrix equation, $U^{2}=I$.
We can again rearrange the expression for the $i,j$ element of $I$
in this expression ($i<j$)
\[
\sum_{k=i}^{j}u_{ik}u_{kj}=0
\]
so as to isolate 
\begin{equation}
u_{ij}=-\frac{\sum_{k=i+1}^{j-1}u_{ik}u_{kj}}{u_{ii}+u_{jj}}\;.\label{eq:element-higham-parlett}
\end{equation}
If $u_{ii}$ and $u_{jj}$ have the opposite sign (a $1$ and a $-1$),
this expression breaks down. However, in this case the signs of $t_{ii}$
and $t_{jj}$ are also different, so the plain Parlett recurrence
(Equation\ref{alg:parlett}) can be safely used. When both $u_{ii}+u_{jj}\neq0$
and $t_{ii}-t_{jj}\neq0$ we prefer to compute $u_{ij}$ using Equation~\ref{alg:parlett-sylvester}
rather than using Equation~\ref{alg:parlett} becasue $|u_{ii}+u_{jj}|=2$
whereas $|t_{ii}-t_{jj}|$ can be small (even if both $t_{ii}$ and
$t_{jj}$ are far from zero). Algorithm~\ref{alg:parlett-higham}
shows the details of this approach.

\begin{algorithm}
\protect\caption{\label{alg:parlett-higham}The Parlett-Higham substitition algorithm
for the matrix sign function.}

\begin{algorithmic}[1]

\STATE$\text{Compute a (complex) Schur decomposition }A=QTQ^{*}$.

\STATE\textbf{for }$i=1:n,\; u_{ii}=\text{sign}(t_{ii})$

\STATE\textbf{for $j=2:n$}

\STATE~~~\textbf{for }$i=j-1:-1:1$

\STATE~~~~~\textbf{if $u_{ii}+u_{jj}=0$}

\STATE~~~~~~~\textbf{then $u_{ij}=t_{ij}\frac{u_{ii}-u_{jj}}{t_{ii}-t_{jj}}+\frac{\sum_{k=i+1}^{j-1}\left(u_{ik}t_{kj}-t_{ik}u_{kj}\right)}{t_{ii}-t_{jj}}$}

\STATE~~~~~~~\textbf{else $u_{ij}=-\frac{\sum_{k=i+1}^{j-1}u_{ik}u_{kj}}{u_{ii}+u_{jj}}$}

\STATE\textbf{~~~end}

\STATE\textbf{end}

\end{algorithmic}
\end{algorithm}

\section{\label{sec:Arithmetic-Efficiency}Arithmetic Efficiency}

Interestingly, the arithmetic efficiency of the two algorithms can
vary considerably (and asymptotically). To design a high-performance
algorithm, we need to choose the most efficient approach for a given
matrix.

The arithmetic complexity of the Parlett-Higham recurrence varies
between $n^{3}/3+o(n^{3})$ and $2n^{3}/3+o(n^{3})$ floating-point
operations (flops). The actual number of operations depends on which
$u_{ij}$s are computed from the equation $U^{2}=I$ (top choice in
Algorithm~\ref{alg:parlett-higham}) and which are computed from
$UT=TU$ (bottom choice), since in the first case the algorithm computes
one inner product on indices ranging from $i+1$ to $j-1$ and in
the second the algorithm computes two such inner products.

The arithmetic complexity of the Parlett-Sylvester algorithm for the
sign function depends on how the eigenvalues of $T$ are initially
ordered along its diagonal. The Schur reordering step (the Bai-Demmel
algorithm or its partitioned variant by Kressner~\cite{KressnerBlockSchurReorderingTOMS})
moves eigenvalues along the diagonal of a triangular matrix by swapping
adjacent eigenvalues using Givens rotations. The number $k$ of swaps
required to group together positive and negative eigenvalues varies
between $0$ and $\frac{n^{2}}{4}$~\cite{Ng:1984:CCM}. The Schur
reordering algorithm performs $12nk$ operations (ignoring low-order
terms)~\cite[Section~7.6.2]{GolubVanLoan4e}, so the cost of this
step varies between nothing (if the eigenvalues are already grouped
by sign) to $3n^{3}$. The $12nk$ operations include those required
to transform $Q$, the orthonormal matrix of Schur vectors.

Once this algorithm reorders the Schur form, it needs to solve a Sylvester
equation for an $n_{-}$by $n_{+}$ off diagonal block, where $n_{-}$
and $n_{+}$ are the numbers of negative and positive eigenvalues.
The number of arithmetic operations required to solve such a Sylvester
equations is 
\[
n^{2}-n\leq n_{-}n_{+}(n_{-}+n_{+})\leq n^{3}/4
\]
(it is easy to see that the extreme cases are $n_{-}=1$ and $n_{-}=n/2$).
We ignore in this analysis the trivial case where all the eigenvalues
are positive or negative, in which the sign is $I$ or $-I$. As in
the first step, the algorithm tends to get more expensive when the
numbers of positive and negative eigenvalues are roughly balanced.

Finally, the algorithm needs to transform the sign of the reordered
matrix to the sign of the input matrix. If this is done by applying
the Givens rotations again, the cost depends on the number of swaps
what the reordering step used. In the best case we need not transform
at all, and in the worst case the cost is cubic.

The critical observation is that in easy cases that require few or
no swaps to reorder the Schur form, the Parlett-Sylvester approach
performs only a quadratic number of floating point operations, whereas
in the worst case, it performs more than $3n^{3}$operations. This
means that this approach can be much more efficient than the Parlett-Higham
approach (if the former performs a quadratic number of operations
and the latter a cubic number) or up to $9$ times less efficient.
Operation counts are not the only determinants of running time, so
the actual performance differences may not be as dramatic, but operation
counts do matter. We address another determinants of performance next.

\section{\label{sec:Communication-Lower-Bounds}Communication Lower Bounds}

We next obtain a communication cost lower bound for Algorithm~\ref{alg:parlett-higham}.
The bound is an application of~ \cite{ballard2011minimizing}, which
extends a technique developed to bound communication in matrix multiplication~\cite{irony2004communication}
to many other computations in linear algebra. The technique embeds
the iteration space of three-nested loops computations into a three
dimension cube and utilizes the Loomis-Whitney~\cite{loomis1949}
inequality to relate operation counts (the volume that the iterations
fill in the cube) to communication requirements (the projections of
the iterations on the input and output matrices).

The lower bound is derived from the computations performed in the
inner loop, lines~5--7. It ignores the computations in line~2 (which
can only increase the total communication cost). Note that either
half or more of the executions of line 5 take the ``then'' branch
(line~6) or half or more take the ``else'' branch on line 7.

We analyze first the second case, in which at least half the time
we have $u_{ii}+u_{jj}\neq0$. We map the computation in line~7 to
Equation 2.1 in~\cite{ballard2011minimizing}. In particular, we
map $u_{ik}$ here to $a(i,k)$ there, $u_{kj}$ to $b(k,j)$, and
$u_{ij}$ to $c(i,j)$. We map the scalar multiplication of $u_{ik}$
by $u_{kj}$ to the abstract function $g_{i,j,k}(\cdot,\cdot)$ in
\cite[Equation~2.1]{ballard2011minimizing}, and the summation and
scaling of the sum by $(u_{ii}+u_{jj})^{-1}$ to the abstract function
$f_{i,j}$. We note that all computed $u_{ij}$ are part of the algorithm's
output, so none of them is discarded; this implies, in the terminology
of~\cite{ballard2011minimizing}, that there are no $R2/D2$ intermediate
results. By applying Theorem 2.2 of~\cite{ballard2011minimizing},
we have,
\begin{cor}
\label{cor:G1}Let $G_{1}$ be the number of arithmetic operations
computed in line 7 of Algorithm~\ref{alg:parlett-higham}. Let $M$
be the cache size. Then the communication cost (number of words transferred
between the cache and main memory) in algorithm~\ref{alg:parlett-higham}
is at least $G_{1}/(8\sqrt{M})-M$.
\end{cor}
We now analyze the communication required to perform the operations
in line~6 of the algorithm, when $u_{ii}+u_{jj}=0$. We again apply
Equation~2.1 and Theorem~2.2 of~\cite{ballard2011minimizing}.
Let $a(i,k)$ there be our $u_{i,k}$, let $b(k,j)$ there be our
$t_{k,j}$, and let $c(i,j)$ there be our $u_{i,j}$. Further, let
$g_{i,j,k}(\cdot,\cdot)$ function be scalar multiplication $u_{ik}\cdot t_{k,j}$,
and $f_{i,j}$ function be the computation of $u_{ij}$, which calls
to $g_{ijk}$. Again we note that all computed $u_{ij}$ are part
of the algorithm's output, so none of them is discarded. We also note
that we can impose writes on the the $n^{2}$ elements of $T$, (see
Section 3.4 of~\cite{ballard2011minimizing}), loosing at most $\Theta(n^{2})$
of the lower bound. Thus, using the terminology of~\cite{ballard2011minimizing},
there are no $R2/D2$ arguments.

By applying Theorem~2.2 of~\cite{ballard2011minimizing}, we have,
\begin{cor}
\label{cor:G2}Let $G_{2}$ be the number of arithmetic operations
perfromed in line 6 of Algorithm~\ref{alg:parlett-higham}. Let $M$
be the cache size. Then the communication cost of the algorithm is
at least $\Omega(G/\sqrt{M}-M-\Theta(n^{2}))$.
\end{cor}
Let $G$ be the total number of arithmetic operations performed in
the doubly-nested loop of Algorithm~\ref{alg:parlett-higham}. Recall
that $\max\{G_{1},G_{2}\}\geq G/2$. Combining Corollary~\ref{cor:G1}
and Corollary~\ref{cor:G2}, we conclude that 
\begin{thm}
Let $G=\Theta(n^{3})$ be the number of arithmetic operations computed
in lines 5--7 of Algorithm~\ref{alg:parlett-higham}, and let $M$
be the size of the cache. The communication cost of Algorithm~\ref{alg:parlett-higham}
is $\Omega(G/\sqrt{M}-M-\Theta(n^{2}))$. Assuming $M<n^{2}$, the
cost is $\Omega(G/\sqrt{M})=\Omega(n^{3}/\sqrt{M})$.
\end{thm}

\section{\label{sec:Communication-Efficient-Algs}Communication-Efficient
Algorithms}

We now propose communication efficient variants of both algorithmic
approaches. We begin with the Parlett-Sylvester approach, which is
more straightforward.

\subsection{Communication-Efficient Parlett-Sylvester Solver}

This approach calls two subroutines: a Schur reordering subroutine
and a Sylverster-equation solver. Fortunatly, communication efficient
variants of both algorithms have been developed. Kressner~\cite{KressnerBlockSchurReorderingTOMS}
developed a communication-efficient variant of the Bai-Demmel reordering
algorithm. Jonsson and Kågström~\cite{Jonsson:2002:RBAb} developped
RECSY, a recursive communication-efficient Sylvester solver.

We have implemented this algorithmic approach in two ways. One calls
xTRSEN, LAPACK's implementation of the Bai and Demmel algorithm that
operates on rows and columns and ignores communication efficiency,
and xTRSYL, LAPACK's Sylvester solver, similarly not communication
efficient. The other calls communication-efficient codes by Kressner
and by Jonsson and Kågström. We use the first LAPACK-based implementation
to evaluate the performance improvement achieved by the new communication-efficient
approach.

\subsection{Communication Efficient Parlett-Higham Solvers}

The communication-efficient algorithm is a recursion that is based
on a nested partitioning of the index set $\{1,2,\ldots,n\}$. The
recursion is somewhat more complex than the recursion for simpler
matrix algorithms (e.g., Cholesky). To present it and to prove its
correctness, we introduce a notation for the nested partitioning and
for sums over subsets of a partition.
\begin{defn}
A \emph{nested partitioning} of $\{1,2,\ldots,n\}$ is a collection
of index sets $p=\{P^{(0)},P^{(1)},\ldots P^{(L)}\}$ such that $P^{(0)}=\{1\}$
and if $P^{(\ell)}=\{i_{1},i_{2}\ldots,i_{m}\}$ then $i_{1}<i_{2}<\cdots<i_{m}$
and $P^{(\ell-1)}=\{i_{1},i_{3},i_{5},\ldots i_{m}\}$ or $P^{(\ell-1)}=\{i_{1},i_{3},i_{5},\ldots i_{m-1}\}$.
\end{defn}
Note that the definition implies that $i_{1}=1$. The indexes in a
partition represent the beginnings of a block of row/column indexes.
For example, $P^{(\ell)}=\{i_{1},i_{2}\ldots,i_{m}\}$ represent the
partitioning of the range $1\colon n$ (in Matlab notation) into $i_{1}\colon i_{2}-1=1\colon i_{2}-1$,
$i_{2}\colon i_{3}-1$, etc. 

For example, let $n=1000$ and let
\begin{eqnarray*}
P^{(0)} & = & \{1\}\\
P^{(1)} & = & \{1,500\}\\
P^{(2)} & = & \{1,250,500,750\}\\
P^{(3)} & = & \{1,125,250,375,500,625,750,875\}\;.
\end{eqnarray*}
We use nested partitions to denote blocks of vectors and matrices.
Using the example above, we can denote blocks of a vector $v$ and
a matrix $A$ by
\begin{eqnarray*}
v_{250}^{(3)} & = & v_{250\colon374}\\
v_{250}^{(2)} & = & v_{250\colon499}\\
A_{250,625}^{(3)} & = & A_{250\colon499,625:749}
\end{eqnarray*}
and so on. In this notation, a block of indices at level $\ell$ must
start at some $i_{j}\in P^{(\ell)}$, and it ends at $i_{j+1}-1$.
We now define a function that allows us to iterate over ranges in
a given partition. 
\begin{defn}
Let $P$ be a nested partitioning and let $P^{(\ell)}=\{i_{1},i_{2}\ldots,i_{m}\}$.
The function $\eta:P^{(\ell)}\rightarrow P^{(\ell)}\cup\{n+1\}$ returns
the start index of the next range in a given partition 
\[
\eta^{(\ell)}(i_{j})=i_{j+1}\text{ (in }P^{(\ell)}\text{)\;.}
\]
For completeness, we define 
\[
\eta^{(\ell)}(i_{m})=n+1\;,
\]
so that subtracting $1$ from the next range always gives the last
element in the current range. We also define the function $\pi$ that
returns the \emph{previous} range,
\[
\pi^{(\ell)}(i_{j})=i_{j-1}
\]
and 
\[
\pi^{(\ell)}(n+1)=i_{m}\;.
\]

\end{defn}
We can now define how vectors and matrices are partitioned, as well
as sum over ranges in a partition.
\begin{defn}
Let $P$ be a nested partition of $\{1,2,\ldots,n\}$, let $v$ be
an $n$ vector and let $A$ be an $n$-by-$n$ matrix. Let $i,j\in P^{(\ell)}$.
We denoted
\[
v_{i}^{(\ell)}=\begin{bmatrix}v_{i}\\
\vdots\\
v_{\eta^{(\ell)}(i)-1}
\end{bmatrix}
\]
and
\[
A_{i,j}^{(\ell)}=\begin{bmatrix}A_{i,j} & \cdots & A_{i,\eta^{(\ell)}(j)-1}\\
\vdots\\
A_{\eta^{(\ell)}(i)-1,j} & \cdots & A_{\eta^{(\ell)}(i)-1,\eta^{(\ell)}(j)-1}
\end{bmatrix}\;.
\]

\end{defn}
Clearly, $v_{i}^{(\ell)}=\begin{bmatrix}v_{i}^{(\ell+1)} & v_{\eta^{(\ell+1)}(i)}^{(\ell+1)}\end{bmatrix}^{T}$
and similarly for matrices. We also need the reverse notation. That
is, we abuse the notation mildly and denote
\begin{eqnarray*}
\left(v_{i}^{(\ell)}\right)_{i}^{(\ell+1)} & = & v_{i}^{(\ell+1)}\\
\left(v_{i}^{(\ell)}\right)_{\eta^{(\ell+1)}(i)}^{(\ell+1)} & = & v_{\eta^{(\ell+1)}(i)}^{(\ell+1)}
\end{eqnarray*}
and similarly for matrices.
\begin{defn}
Let $P$ be a nested partitioning and let $P^{(\ell)}=\{i_{1},i_{2}\ldots,i_{m}\}$,
let $s\in P^{(\ell)}$, and let $e\in P^{(\ell)}$ or $e=n+1$. We
define 
\[
\sum_{j=s}^{e-1}v_{j}^{(\ell)}=\begin{cases}
0 & s>e\\
\sum_{j=s}^{\eta^{(\ell)}(s)-1}v_{j}+\sum_{j=\eta^{(\ell)}(s)}^{e-1}v_{j}^{(\ell)} & \text{otherwise.}
\end{cases}
\]
The sum consists of all the elements of $v$ starting at the beginning
of a range in $P^{(k)}$ and ending just before another range in $P^{(k)}$
starts. Note that the first sum on the right hand side is a sum over
scalars that iterates over consecutive integer indexes, whereas the
second sum is defined (recursively) over sums of ranges. The superscript
$(\ell)$ on the argument $v$ (or the lack of superscript) indicates
the type of the sum.
\end{defn}
The following lemma relates sums over ranges in adjacent partitions
in a nest.
\begin{lem}
Let $P$ be a nested partitioning and let $P^{(\ell)}=\{i_{1},i_{2}\ldots,i_{m}\}$,
let $s\in P^{(\ell)}$, and let $e\in P^{(\ell)}$ or $e=n+1$. The
following relation hold,
\[
\sum_{j^{(k)}=s}^{e-1}v_{j^{(k)}}^{(\ell)}=\begin{cases}
\sum_{j=s}^{e-1}v_{j}^{(\ell-1)} & \text{if }s\text{ is odd and }e\text{ is even}\\
v_{s}^{(\ell)}+\sum_{j=\eta^{(\ell)}(s)}^{e-1}v_{j}^{(\ell-1)} & \text{if }s\text{ is even and }e\text{ is even}\\
\sum_{j=s}^{e-1}v_{j}^{(\ell-1)}+v_{\pi^{(\ell)}\left(e\right)}^{(\ell)} & \text{if }s\text{ is odd and }e\text{ is odd}\\
v_{i_{s}}^{(\ell)}+\sum_{j=\eta^{(\ell)}\left(i_{s}\right)}^{\pi^{(k)}\left(i_{e}\right)-1}v_{j}^{(\ell-1)}+v_{\pi^{(\ell)}\left(i_{e}\right)}^{(\ell)} & \text{if }s\text{ is even and }e\text{ is odd.}
\end{cases}
\]

\end{lem}
The Higham-Parlett recurrence is based on the observation that the
sign $U$ of $T$ satisifes both $TU=UT$ and $U^{2}=I$. Neither
of these equations alone defines all the elements of $U$ but together
they do. We partition $U$ and $T$ into block matrices with square
diagonal blocks using a nested partition $P$. The blocks also satisfy
the equations, so for any $k$ in the nest, 
\begin{eqnarray*}
\left(TU\right)_{ij}^{(\ell)} & = & \left(UT\right)_{ij}^{(\ell)}\\
\left(UU\right)_{ij}^{(\ell)} & = & I_{ij}^{(\ell)}\;.
\end{eqnarray*}
which expands into
\begin{eqnarray*}
T_{ii}^{(\ell)}U_{ij}^{(\ell)}-U_{ij}^{(\ell)}T_{jj}^{(\ell)} & = & U_{ii}^{(\ell)}T_{ij}^{(\ell)}-T_{ij}^{(\ell)}U_{jj}^{(\ell)}+\sum_{k=\eta^{(\ell)}(i)}^{j-1}\left(U_{ik}^{(\ell)}T_{kj}^{(\ell)}-T_{ik}^{(\ell)}U_{kj}^{(\ell)}\right)\\
U_{ii}^{(\ell)}U_{ij}^{(\ell)}+U_{ij}^{(\ell)}U_{jj}^{(\ell)} & = & I_{ij}^{(\ell)}-\sum_{k=\eta^{(\ell)}(i)}^{j-1}U_{ik}^{(\ell)}U_{kj}^{(\ell)}\;.
\end{eqnarray*}
We denote the sums on the right by 
\[
X_{ij}^{(\ell)}=\sum_{k=\eta^{(\ell)}(i)}^{j-1}\left(U_{ik}^{(\ell)}T_{kj}^{(\ell)}-T_{ik}^{(\ell)}U_{kj}^{(\ell)}\right)
\]
and
\[
Y_{ij}^{(\ell)}=\sum_{k=\eta^{(\ell)}(i)}^{j-1}U_{ik}^{(\ell)}U_{kj}^{(\ell)}\;.
\]

We now related the blocks of $X$ and $Y$ at level $\ell$ to those
at level $\ell+1$. The easiest one is the $(2,1)$ block,
\begin{eqnarray*}
Y_{\eta^{(\ell+1)}(i),j}^{(\ell+1)} & = & \sum_{k=\eta^{(\ell+1)}(\eta^{(\ell+1)}(i))}^{j-1}U_{ik}^{(\ell+1)}U_{kj}^{(\ell+1)}\\
 & = & \sum_{k=\eta^{(\ell)}(i)}^{j-1}U_{ik}^{(\ell+1)}U_{kj}^{(\ell+1)}\\
 & = & \left(\sum_{k=\eta^{(\ell)}(i)}^{j-1}U_{ik}^{(\ell)}U_{kj}^{(\ell)}\right)_{\eta^{(\ell+1)}(i),j}^{(\ell+1)}\\
 & = & \left(Y_{i,j}^{(\ell)}\right)_{\eta^{(\ell+1)}(i),j}^{(\ell+1)}\;.
\end{eqnarray*}
In the $(2,1)$ and $(2,2)$ blocks, we need to add a contribution
at the $\ell+1$ level,
\begin{eqnarray*}
Y_{i,j}^{(\ell+1)} & = & \sum_{k=\eta^{(\ell+1)}(i)}^{j-1}U_{ik}^{(\ell+1)}U_{kj}^{(\ell+1)}\\
 & = & U_{i,\eta^{(\ell+1)}(i)}^{(\ell+1)}U_{\eta^{(\ell+1)}(i),j}^{(\ell+1)}+\sum_{k=\eta^{(\ell+1)}(\eta^{(\ell+1)}(i))}^{j-1}U_{ik}^{(\ell+1)}U_{kj}^{(\ell+1)}\\
 & = & U_{i,\eta^{(\ell+1)}(i)}^{(\ell+1)}U_{\eta^{(\ell+1)}(i),j}^{(\ell+1)}+\sum_{k=\eta^{(\ell)}(i)}^{j-1}U_{ik}^{(\ell+1)}U_{kj}^{(\ell+1)}\\
 & = & U_{i,\eta^{(\ell+1)}(i)}^{(\ell+1)}U_{\eta^{(\ell+1)}(i),j}^{(\ell+1)}+\left(\sum_{k=\eta^{(\ell)}(i)}^{j-1}U_{ik}^{(\ell)}U_{kj}^{(\ell)}\right)_{i,j}^{(\ell+1)}\\
 & = & U_{i,\eta^{(\ell+1)}(i)}^{(\ell+1)}U_{\eta^{(\ell+1)}(i),j}^{(\ell+1)}+\left(Y_{i,j}^{(\ell)}\right)_{i,j}^{(\ell+1)}\;,
\end{eqnarray*}
and
\[
Y_{\eta^{(\ell+1)}(i),\eta^{(\ell+1)}(j)}^{(\ell+1)}=\left(Y_{i,j}^{(\ell)}\right)_{\eta^{(\ell+1)}(i),\eta^{(\ell+1)}(j)}^{(\ell+1)}+U_{\eta^{(\ell+1)}(i),j}^{(\ell+1)}U_{j,\eta^{(\ell+1)}(j)}^{(\ell+1)}\;.
\]
The $(1,2)$ block requires two contributions from level $\ell+1$,
\[
Y_{i,\eta^{(\ell+1)}(j)}^{(\ell+1)}=U_{i,\eta^{(\ell+1)}(i)}^{(\ell+1)}U_{\eta^{(\ell+1)}(i),\eta^{(\ell+1)}(j)}^{(\ell+1)}+\left(Y_{i,j}^{(\ell)}\right)_{i,\eta^{(\ell+1)}(j)}^{(\ell+1)}+U_{i,j}^{(\ell+1)}U_{j,\eta^{(\ell+1)}(j)}^{(\ell+1)}\;.
\]
The expressions for the blocks of $X$ at level $\ell+1$ are similar.

We can now present the algorithm, which we split into three procedures.
The top-level procedure \emph{sign }allocates $U$, $X$ and $Y$
and zeros $X$ and $Y$. It calls a recursive procedure that computes
a diagonal block of $U$ at level $\ell=0$ called \emph{sign-diagonal}.
Sign-diagonal calls itself recursively to compute the two diagonal
blocks at level $\ell+1$ and a third procedure, \emph{sign-offdiagonal},
which computes an offdiagonal block of $U$. Sign-offdiagonal works
by calling itself four times on the four sub-blocks at the next level. 

\begin{algorithm}
\protect\caption{\label{alg:recursive-toplevel}A procedure that allocates two auxiliary
matrices, $X$ and $Y$, and calls the recursive algorithm to compute
the sign of a triangular matrix $T$.}

\begin{algorithmic}[1]

\STATE$\text{function }U=\text{sign}(T)$

\STATE$\text{allocate }n\text{-by-}n\text{ upper triangular matrices }U$,
$X$ , and $Y$

\STATE$\text{set }X=X^{(0)}=0\text{, }Y=Y^{(0)}=0$

\STATE\textbf{$\text{sign-diagonal}(1,0,T,U,X,Y)$}

\STATE$\text{return }U$

\end{algorithmic}
\end{algorithm}

\begin{algorithm}
\protect\caption{\label{alg:recursive-diagonal}A recursive algorithm to compute a
diagonal block $U_{ii}^{(\ell)}$ of the sign $U$ of $T$. We assume
that the arguments are passed by reference and that the code modifies
elements of arguments $U$, $X$, and $Y$. }

\begin{algorithmic}[1]

\STATE$\text{function }\text{sign-diagonal}(i,\ell,T,U,X,Y)$

\STATE\hspace*{0cm}if $U_{ii}^{(\ell)}$ is $1$-by-$1$ then $U_{ii}^{(\ell)}=u_{ii}=\text{sign}(t_{ii})$.

\STATE$\text{otherwise,}$

\STATE\textbf{$\text{sign-diagonal}(i,\ell+1,T,U,X,Y)$}

\STATE\textbf{$\text{sign-diagonal}(\eta^{(\ell+1)}(i),\ell+1,T,U,X,Y)$}

\STATE\textbf{$\text{sign-offdiagonal}(i,\eta^{(\ell+1)}(i),\ell+1,T,U,X,Y)$}

\STATE$\text{return}$

\end{algorithmic}
\end{algorithm}
The auxiliary algorithm is a bit more complex. 

\begin{algorithm}
\protect\caption{\label{alg:recursive-offdiagonal}A recursive implementation of Parlett-Higham
algorithm to compute an off-diagonal block $U_{ij}^{(\ell)}$. We
again assume that the arguments are passed by reference. Elements
of $U$ that have been computed in previous steps (calls to level
$\ell+1$) are marked in red to emphasize dependencies.}

\begin{algorithmic}[1]

\STATE$\text{function }\text{sign-offdiagonal}(i,j,\ell,T,U,X,Y)$

\STATE\hspace*{0cm}if $U_{ij}^{(\ell)}$ is $1$-by-$1$ then $U_{ij}^{(\ell)}=u_{ij}$,
which we compute as

\STATE\hspace*{0cm}\quad{}$u_{ij}=\begin{cases}
\begin{array}{l}
\frac{-y_{ij}}{u_{ii}+u_{jj}},\\
t_{ij}\frac{u_{ii}-u_{jj}}{t_{ii}-t_{jj}}+\frac{x_{ij}}{t_{ii}-t_{jj}},
\end{array} & \begin{array}{c}
u_{ii}+u_{jj}\neq0\\
u_{ii}+u_{jj}=0
\end{array}\end{cases}$

\STATE$\text{and return. otherwise,}$

\STATE\textbf{$\text{sign-offdiagonal}(\eta^{(\ell+1)}(i),j,\ell+1,T,U,X,Y)$}

\STATE$X_{i,j}^{(\ell+1)}=\left(X_{i,j}^{(\ell)}\right)_{i,j}^{(\ell+1)}+\left(U_{i,\eta^{(\ell+1)}(i)}^{(\ell+1)}{\color{black}{\color{red}{\color{black}T_{\eta^{(\ell+1)}(i),j}^{(\ell+1)}}}}-T_{i,\eta^{(\ell+1)}(i)}^{(\ell+1)}{\color{red}U_{\eta^{(\ell+1)}(i),j}^{(\ell+1)}}\right)$

\STATE$Y_{i,j}^{(\ell+1)}=\left(Y_{i,j}^{(\ell)}\right)_{i,j}^{(\ell+1)}+U_{i,\eta^{(\ell+1)}(i)}^{(\ell+1)}{\color{red}U_{\eta^{(\ell+1)}(i),j}^{(\ell+1)}}$

\STATE\textbf{$\text{sign-offdiagonal}(i,j,\ell+1,T,U,X,Y)$}

\STATE$X_{\eta^{(\ell+1)}(i),\eta^{(\ell+1)}(j)}^{(\ell+1)}=X_{\eta^{(\ell+1)}(i),\eta^{(\ell+1)}(j)}^{(\ell+1)}+\left({\color{red}U_{\eta^{(\ell+1)}(i),j}^{(\ell+1)}}T_{j,\eta^{(\ell+1)}(j)}^{(\ell+1)}-{\color{red}{\normalcolor T_{\eta^{(\ell+1)}(i),j}^{(\ell+1)}}}U_{j,\eta^{(\ell+1)}(j)}^{(\ell+1)}\right)$

\STATE$Y_{\eta^{(\ell+1)}(i),\eta^{(\ell+1)}(j)}^{(\ell+1)}=\left(Y_{i,j}^{(\ell)}\right)_{\eta^{(\ell+1)}(i),\eta^{(\ell+1)}(j)}^{(\ell+1)}+{\color{red}U_{\eta^{(\ell+1)}(i),j}^{(\ell+1)}}U_{j,\eta^{(\ell+1)}(j)}^{(\ell+1)}$

\STATE\textbf{$\text{sign-offdiagonal}(\eta^{(\ell+1)}(i),\eta^{(\ell+1)}(j),\ell+1,T,U,X,Y)$}

\STATE$X_{\eta^{(\ell+1)}(i),\eta^{(\ell+1)}(j)}^{(\ell+1)}=X_{\eta^{(\ell+1)}(i),\eta^{(\ell+1)}(j)}^{(\ell+1)}+\left(U_{i,\eta^{(\ell+1)}(i)}^{(\ell+1)}{\color{red}{\normalcolor T_{\eta^{(\ell+1)}(i),\eta^{(\ell+1)}(j)}^{(\ell+1)}}}-T_{i,\eta^{(\ell+1)}(i)}^{(\ell+1)}{\color{red}U_{\eta^{(\ell+1)}(i),\eta^{(\ell+1)}(j)}^{(\ell+1)}}\right)+\left({\color{red}U_{i,j}^{(\ell+1)}}T_{j,\eta^{(\ell+1)}(j)}^{(\ell+1)}-{\color{red}{\normalcolor T_{i,j}^{(\ell+1)}}}U_{j,\eta^{(\ell+1)}(j)}^{(\ell+1)}\right)$

\STATE$Y_{i,\eta^{(\ell+1)}(j)}^{(\ell+1)}=\left(Y_{i,j}^{(\ell)}\right)_{i,\eta^{(\ell+1)}(j)}^{(\ell+1)}+U_{i,\eta^{(\ell+1)}(i)}^{(\ell+1)}{\color{red}U_{\eta^{(\ell+1)}(i),\eta^{(\ell+1)}(j)}^{(\ell+1)}}+{\color{red}U_{i,j}^{(\ell+1)}}U_{j,\eta^{(\ell+1)}(j)}^{(\ell+1)}$

\STATE\textbf{$\text{sign-offdiagonal}(i,\eta^{(\ell+1)}(j),\ell+1,T,U,X,Y)$}

\STATE$\text{return}$

\end{algorithmic}
\end{algorithm}

\subsection{\label{sub:arithmetic-efficient-recursive-higham}Improving The Arithmetic
Complexity}

Algorithm~\ref{alg:recursive-toplevel} performs $n^{3}$ arithmetic
operations, more than the $n^{3}/3$ to $2n^{3}/3$ operations that
the Parlett-Higham recurrence performs. This happens because extended-sign
computes both 
\[
\sum_{k=i+1}^{j-1}u_{ik}u_{kj}\text{ and }\sum_{k=i+1}^{j-1}\left(u_{ik}t_{kj}-t_{ik}u_{kj}\right)
\]
for every $i<j$, whereas Algorithm~\ref{alg:parlett-higham} only
computes one of the two for a particular $i,j$. In other words, the
algorithm computes all the entries of both $X$ and $Y$ but it does
not actually use all of them later. For a given position $i,j$, only
one of $x_{ij}$ and $y_{ij}$ is needed, the one that the sign-offdiagonal
function needs. If $\bar{u}_{ii}+\bar{u}_{jj}=0$, we need $x_{ij}$;
otherwise, it is $y_{ij}$.

We can improve the arithmetic complexity of the algorithm by computing
only one of $x_{ij}$ and $y_{ij}$. More specifically, when calculating
the contributions to $X$ and $Y$ in between the recursive calls
in sign-offdiagonal, we only compute elements of the $X$ argument
that are actually needed and only elements of $Y$ that are actually
needed. In practice, we can only keep one matrix $Z$ and decide on
the method of calculating $z_{ij}$ based on the values $u_{ii}$
and $u_{jj}$.

This approach performs fewer arithmetic operations (by a factor of
$2$ to $3$), but it prevents us from using existing matrix multiplication
codes (e.g., xGEMM), so it is unlikely to be fast in practice. We
have implemented this algorithm but the experiments below demonstrate
that it is indeed slow.

\section{\label{sec:Experimental-Results}Experimental Results}

We evaluated several different algorithms experimentally. We implemented
the algorithms in C and called them from Matlab for testing and we
used the BLAS and LAPACK libraries that are bundled with Matlab. We
used Matlab R2013A which uses Intel's Math Kernel Library Version
10.3.11 for the BLAS and LAPACK and is based on LAPACK version 3.4.1.

We conducted the experiments on a quad-core desktop computer running
Linux. The computer had 16~GB of RAM and an Intel~i7-4770 CPU processor
running at 3.40~GHz. Some of the experiments used only one core (using
\texttt{maxNumCompThreads(1)} in Matlab) and some used all four (same
function with argument 4), but only in BLAS routines. Runs that used
4 cores are labeled \emph{MT }in the graphs below.

We tested all the algorithms on random triangular matrices with a
prescribed inertia. We generated the matrices by creating random real
square matrices with elements that are distributed uniformly in $[-50,50]$,
computing their complex Schur form, and taking the real part of the
Schur form. This generates matrices with roughly balanced inertia.
In the experiments reported below, the fraction of negative eigenvalues
ranged from 48\% to 54\% on the smallest matrices (dimension 50),
from 49\% to 51\% on the next smallest dimension (657), and even narrower
on larger matrices. In some of the experiments we forced the number
negative eigenvalues to a prescribed number $k$. We did this by keeping
the absolute values of the diagonal elements of the random triangular
matrix, but forcing their sign to positive in all but a random $k$
positions.

We tested the following algorithms:
\begin{itemize}
\item The Parlett-Higham algorithm (Algorithm~\ref{eq:element-higham-parlett}).
We refer to this algorithm as \emph{Higham} in the graphs below.
\item Two implementations of the Parlett-Sylvester algorithm (specialized
to the sign function). The first implementation uses LAPACK's built-in
routines for reordering the Schur form and for solving the Sylvester
equations. Neither routine is blocked in LAPACK 3.4.1. We refer to
this implementation as \emph{LAPACK Sylvester}.
\item The second implementation of the Parlett-Sylvester algorithm used
RECSY, a recursive Sylvester solver by Jonsson and Kågström~\cite{Jonsson:2002:RBAb},
as well asd a blocked Schur reordering code by Kressner~\cite{KressnerBlockSchurReorderingTOMS}.
\item Our recursive implementation of the Parlett-Higham algorithm (Algorithms~\ref{alg:recursive-toplevel},
\ref{alg:recursive-diagonal}, and \ref{alg:recursive-offdiagonal}).
This implementation calles the BLAS to multiply blocks. Recursion
was used only on blocks with dimension larger than $16$; smaller
diagonal blocks were processed by our element-by-element Parlett-Higham
implementation. We refer to this implementation as \emph{Recursive
Higham} \emph{MM}.
\item A recursive implementation of the arithmetic-efficient Parlett-Higham
algorithm described in Section~\ref{sub:arithmetic-efficient-recursive-higham}.
This implementation does not use the BLAS (as its operations do not
reduce to matrix multiplications). We refer to it as \emph{Recursive
Higham}.
\end{itemize}
\begin{figure}
\begin{centering}
\includegraphics[width=1\columnwidth]{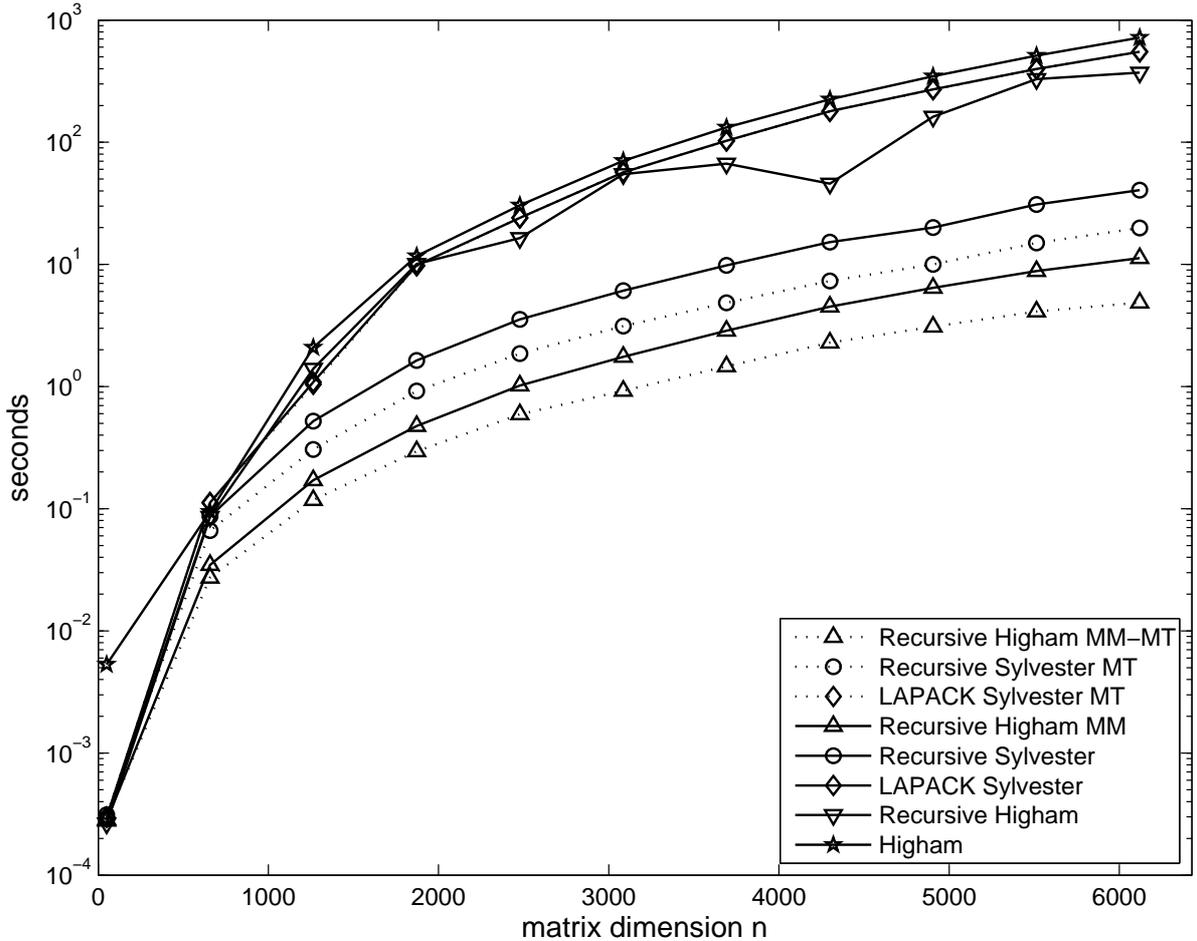}
\par\end{centering}

\protect\caption{\label{fig:time-balanced}Running times on matrices with roughly balanced
inertia.}
\end{figure}

The running times on matrices with roughly balanced inertia are shown
in Figure~\ref{fig:time-balanced}. Our recursive algorithm is the
fastest one, both with and without multithreaded BLAS. The next-best
algorithm is the recursive Patlett-Sylvester algorithm. Like our recursive
algorithm, it uses the BLAS extensively so it benefits from multithreading.
Our recursive but arithmetic-efficient algorithm is fairly slow, because
it does not use the BLAS. The slowest algorithms are the Parlett-Sylvester
implementation that uses LAPACK for Schur reordering and for solving
Sylvester equations and the element-by-element Parlett-Higham algorithm.

\begin{figure}
\begin{centering}
\includegraphics[width=1\columnwidth]{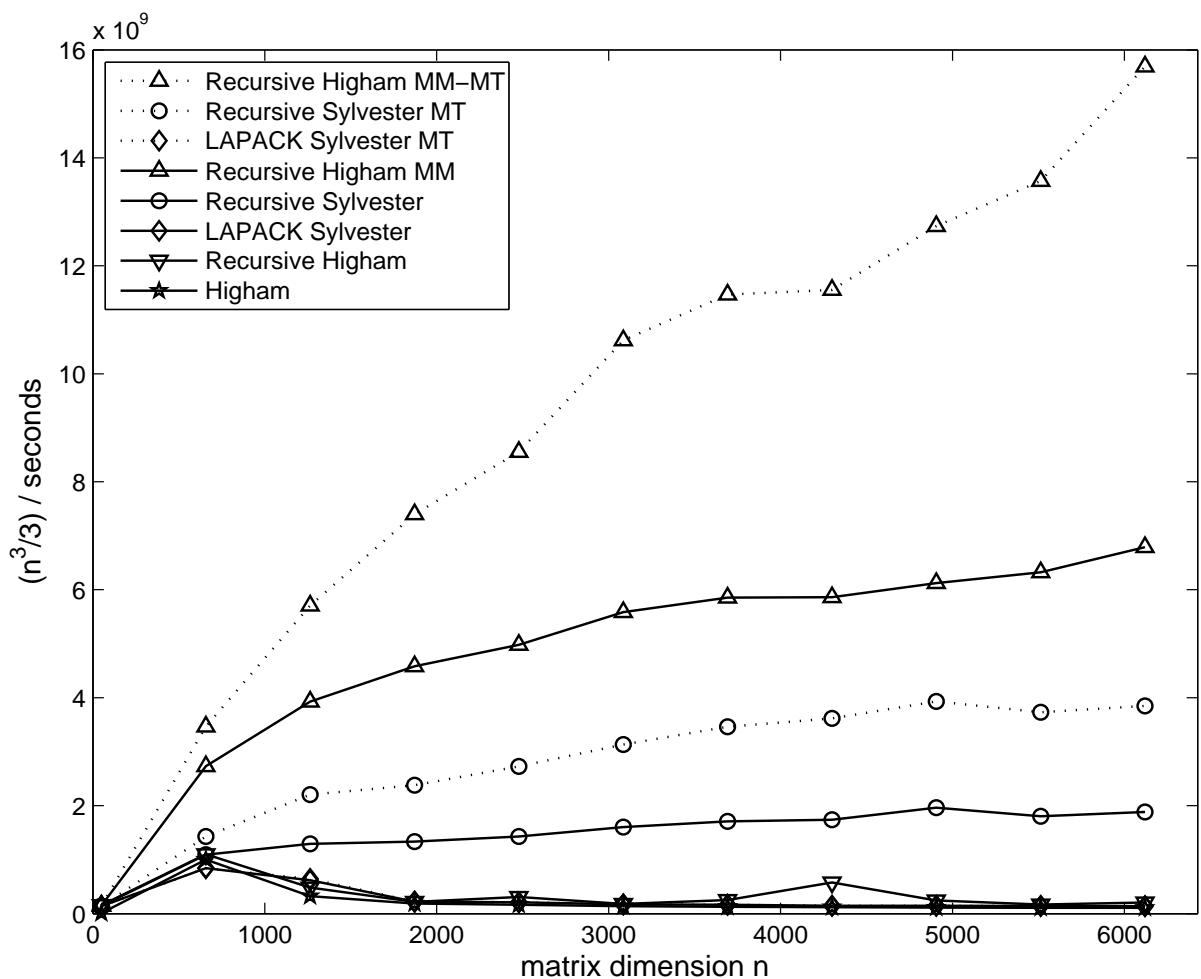}
\par\end{centering}

\protect\caption{\label{fig:rate-balanced}Normalized computational rates on matrices
with roughly balanced inertia. The number $n^{3}/3$ is used for normalization
because the number of arithmetic operations in Higham's algorithm
is between $n^{3}/3+o(n^{3})$ and $2n^{3}/3+o(n^{3})$; the number
of operations in some of the other algorithms is different.}
\end{figure}

Figure~\ref{fig:rate-balanced} puts the same results in a somewhat
more familiar quantitative context. By measuring performance in terms
of normalized floating-point arithmetic rates, the performance of
the algorithms can be directly compared to the performance of other
algorithms (e.g., matrix multiplication) on the same computer. The
rates are normalized relative to $n^{3}/3$ because the number of
arithmetic operations in Higham's algorithm is between $n^{3}/3+o(n^{3})$
and $2n^{3}/3+o(n^{3})$; other algorithms may perform more or less
arithmetic.

Our recursive algorithm always performs $2n^{3}/3+o(n^{3})$; on large
matrices it runs single threaded at a rate of about 12Gflop/s (not
normalized). Multithreading on the quad-core computer speeds up the
algorithm by more than a factor of 2 on large matrices (the speedup
is around 2 rather than 4 because only matrix multiplications exploit
more than one core). The recursive Parlett-Sylvester is about 3 times
slower. The performance of the non-recursive algorithms (and of our
recursive algorithm that does not use the BLAS) is quite dismal.

\begin{figure}
\begin{centering}
\includegraphics[width=1\columnwidth]{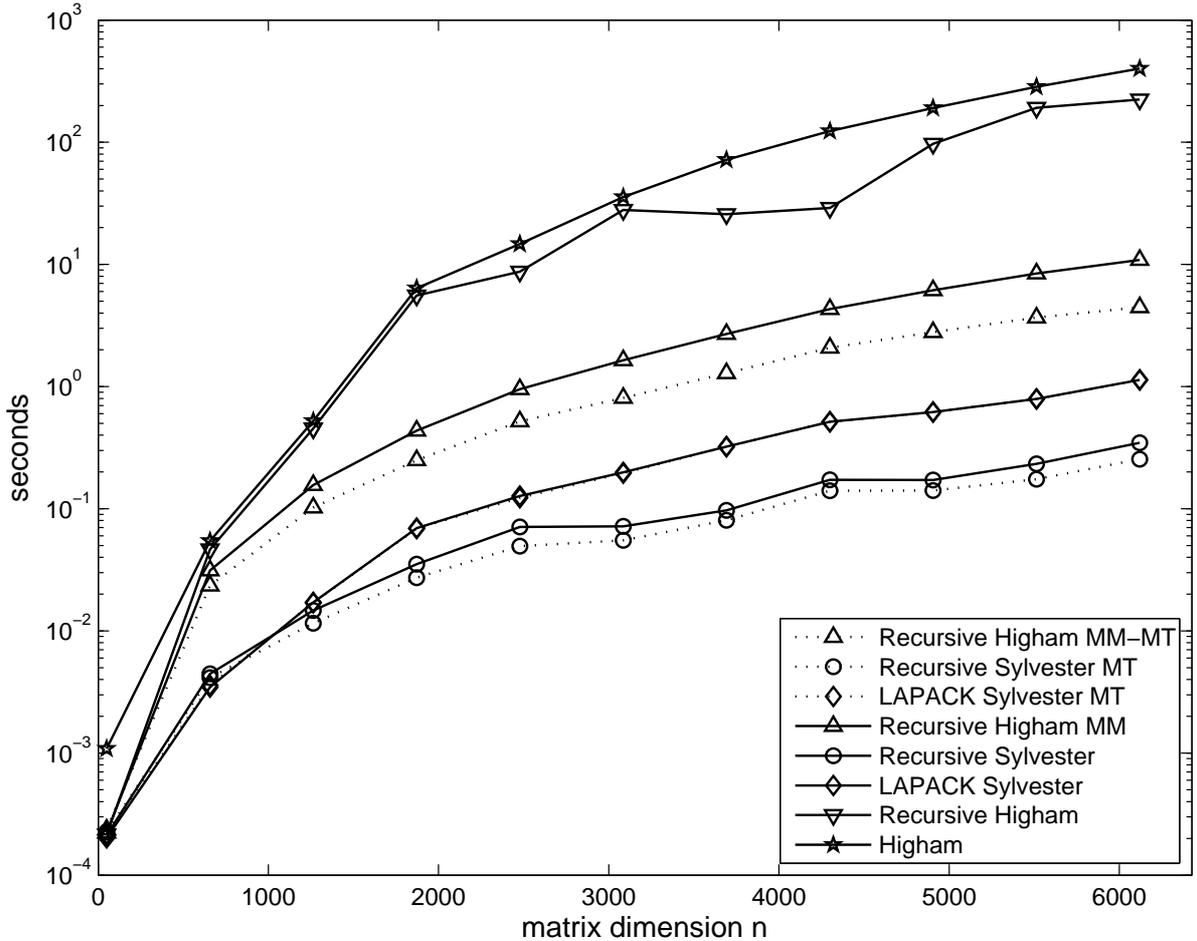}
\par\end{centering}

\protect\caption{\label{fig:time-3}Running times on matrices with Exactly $3$ negative
eigenvalues and $n-3$ positive eigenvalues.}
\end{figure}

\begin{figure}
\begin{centering}
\includegraphics[width=1\columnwidth]{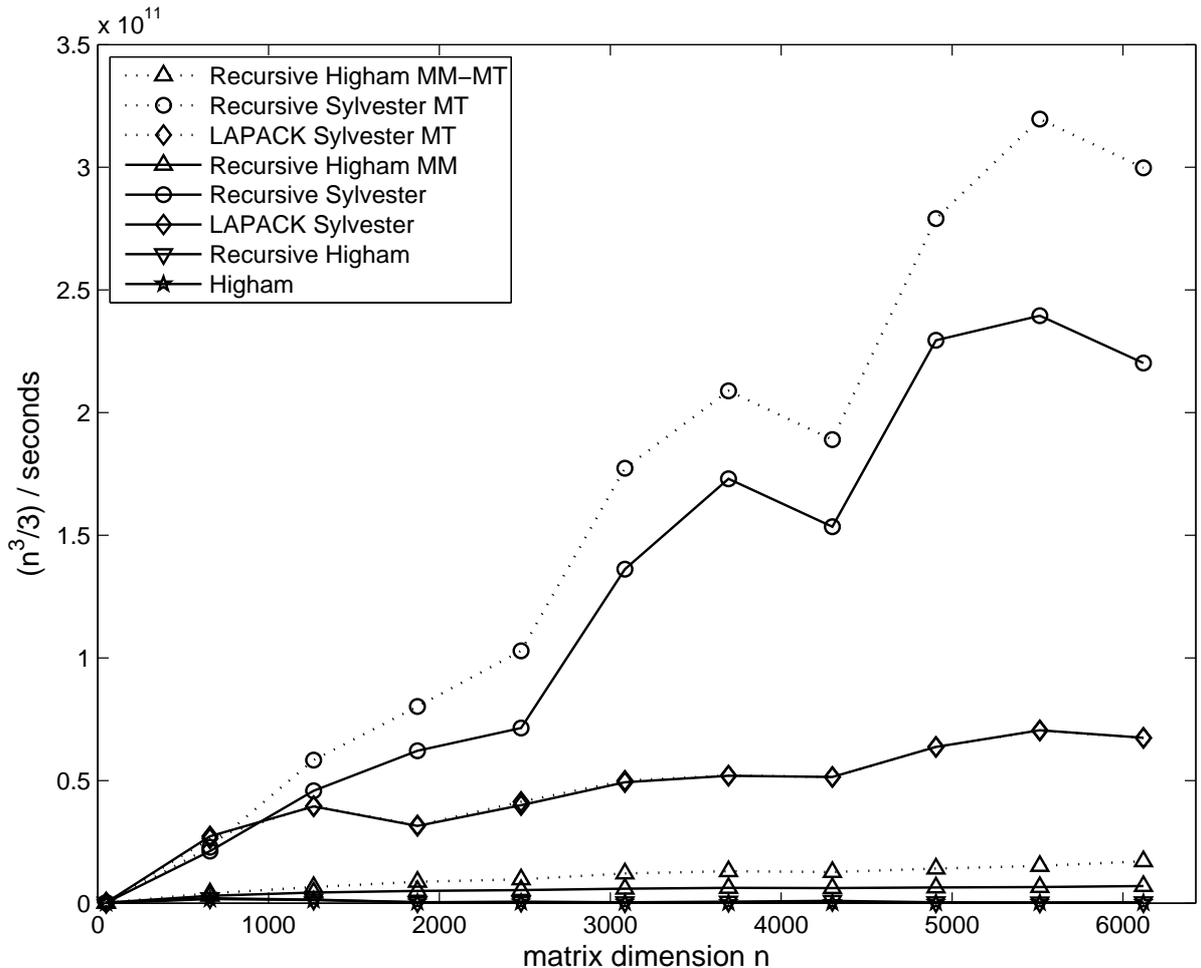}
\par\end{centering}

\protect\caption{\label{fig:rate-3}Normalized computational rates on matrices with
Exactly $3$ negative eigenvalues and $n-3$ positive eigenvalues.}
\end{figure}

When inertia is highly imbalanced, the picture changes. Figure~\ref{fig:time-3}
shows that Parlett-Sylvester algorithms are the fastest on such matrices.
This makes sense, as they only perform $\Theta(n^{2})$ operations,
not $\Theta(n^{3})$ like all the other algorithms. The differences
are quite dramatic. The best single-threaded Parlett-Sylvester algorithm
(\emph{Recursive Sylvester}) runs in 0.35s on matrices of dimension
6120, whereas the fastest single-threaded recursive Higham algorithm
takes 10.9s (more than 30 times slower). Figure~\ref{fig:rate-3}
shows the corresponding normalized rates for completeness, but they
are less interesting because the normalization factor is way off for
the Parlett-Sylvester algorithms.

\section{\label{sec:Conclusions}Conclusions}

The reader may have been somewhat surprised by some aspects of this
work. They also surprised us.

The first surprise is that arithmetic performance (number of operations)
can differ so dramatically between the Parlett-Higham recurrence and
its block variant that we refer to as Parlett-Sylvester. The striking
efficiency of the Parlett-Sylvester approach on matrices with highly
imbalanced inertia is the result of three contributing factors: (1)
the performance of the Schur reordering algorithm depends strongly
on the number of eigenvalue swaps required to order the matrix, (2)
solving Sylverster equations on high-aspect ratio matrices is very
inexpensive, and (3) computing the diagonal blocks in the Parlett-Sylvester
algorithm \emph{for the sign function} is trivial.

This finding implies that a production code for the sign function
should choose between these two algorithms, ideally through an auto-tuning
and/or performance-prediction framework, possibly based on inertia
estimation.

Second, was the difficulty of expressing cleanly the recursive variant
of the Higham-Parlett algorithm. We have tried a number of approaches
based on conventional notational schemes and failed. We resorted to
develop the somewhat complex notation that we present and use in Section~\ref{sec:Communication-Efficient-Algs};
it may seem overly complex, but we found it impossible to present
the algorithm without it. 

The third (and relatively minor) surprise is the benefit of performing
more arithmetic in order to use matrix-multiplication. The arithmetic-efficient
variant of the recursive Parlett-Higham algorithm (Section~\ref{sub:arithmetic-efficient-recursive-higham})
is slower in practice, although it is cache efficient. Rather than
using existing matrix-multiplication routines (xGEMM), it uses a custom
kernel with a condition in the next-to-inner loop. This demonstrated
the performance penalty for trying to do less arithmetic in an algorithm
using a conditional, thus making performance optimization difficult.

\section*{Acknowlegements}

This research was supported in part by grants 863/15, 1878/14, and
1901/14 from the Israel Science Foundation (founded by the Israel
Academy of Sciences and Humanities) , grant 3-10891 from the Ministry
of Science and Technology, Israel. Research is also supported by the
Einstein Foundation and the Minerva Foundation. This paper is supported
by the Intel Collaborative Research Institute for Computational Intelligence
(ICRI-CI). This research was supported by a grant from the United
States-Israel Binational Science Foundation (BSF), Jerusalem, Israel.
This work was supported by the HUJI Cyber Security Research Center
in conjunction with the Israel National Cyber Bureau in the prime
minister's office.

\bibliographystyle{plain}
\bibliography{signfunction,functions-of-matrices}

\end{document}